\documentclass[12pt,
reqno]{amsart}

\usepackage[cp1251]{inputenc}
\usepackage{amsmath,amsxtra,amssymb,eufrak}%{amscd}
\usepackage[all]{xy}
\usepackage[active]{srcltx} %SRC Specials for DVI Searching
\usepackage{mathrsfs}
\usepackage[english]{babel}

\DeclareMathOperator{\LinC}{Lin_\mathbb{C}}

\theoremstyle{plain}
\newtheorem{thm}{Theorem}
\newtheorem{lm}[thm]{Lemma}
\newtheorem{co}[thm]{Corollary}
\newtheorem{pr}[thm]{Proposition}

\theoremstyle{definition}
\newtheorem{rem}[thm]{Remark}

\title{The structure of the linearizer of a connected complex Lie group}

\author{O. Yu. Aristov}
\email{aristovoyu@inbox.ru}

\begin{document}
 \maketitle
 \markright{The structure of the linearizer}

\begin{abstract}
The Morimoto theorem states that each connected abelian complex Lie group $A$ can be
decomposed into the direct product of a group on which all holomorphic functions are constant, finitely
many copies of $\mathbb{C}^\times$ and a vector group. We prove that if~$A$ is the complex linearizer of a connected
complex Lie group then the last factor of the product is trivial.
\end{abstract}

Recall that the \emph{linearizer} of a complex Lie group $G$ is the intersection of the kernels of all holomorphic
finite-dimensional representations. We denote the linearizer of $G$ by $\LinC(G)$. This notion has been
known circa the 1950s, but the main attention of researchers was paid to the conditions for the \emph{linearity} of $G$, i.e., to the triviality of $\LinC(G)$.   At the same time the structure of the linearizer has been not
considered earlier, to the author's knowledge. The aim of the note is to fill this gap. Addressing the
problem, the author pursues his goals, since the description of the linearizer is used in questions about
the homological properties of the algebras of analytical functionals \cite{Ar21}. However, it seems that the result
of this article is of interest in its own right.

The kernel of the adjoint representation of a connected Lie group $G$ coincides with the center of~$G$.
Therefore, the subgroup $\LinC(G)$ is central and, in particular, abelian. Furthermore, $\LinC(G)$  is connected (see Proposition~\ref{Linconn} below). The structure of connected abelian complex Lie groups was studied by
Matsushima and Morimoto in \cite{MM60,Mo65,Mo66}.  Namely, such a group has the following form (here $\mathbb{C}^\times$ is the group of the invertible elements of~$\mathbb{C}$):
\begin{equation}\label{Mdec}
M\times (\mathbb{C}^\times)^k \times\mathbb{C}^m \qquad (k,m\in\mathbb{Z}_+),
\end{equation}
where $M$ is an abelian complex Lie group such that all holomorphic functions on $M$ are constant
\cite[Theorem 3.2]{Mo66} (see \cite{HiNe}, Lemma 15.2.12(ii), Proposition 15.3.4 and Corollary 15.3.5 as well). Note that $M$ need not be compact (see \cite{Mo65} or \cite[Example 15.3.10]{HiNe}).

By \cite[Theorem 1]{Mo65}, if $G$ is a connected complex Lie group, then there is a smallest closed normal
subgroup $M$ such that $G/M$ is a Stein group. In this case all holomorphic functions on $M$ are constant. We refer to $M$ as the \emph{Morimoto subgroup} of~$G$.  In particular, if $G$ is abelian then $M$ in~\eqref{Mdec} is the
Morimoto subgroup of $G$ \cite[Theorem 3.2]{Mo66}.

Our main result states that the last factor is trivial in the decomposition~\eqref{Mdec}  of $\LinC(G)$  More
accurately, the following theorem holds.
 \begin{thm}\label{linarizc}
Let $G$ be a connected complex Lie group, and let $M$ be the Morimoto subgroup of~$G$.
Then $\LinC(G)$ is isomorphic to $M\times(\mathbb{C}^\times)^k$ for some $k\in\mathbb{Z}_+$.
\end{thm}

From this we immediately obtain the linearizer of a Stein group (by the triviality of the Morimoto
subgroup in this case).

\begin{co}\label{linarizcSt}
Let $G$ be a connected Stein group. Then $\LinC(G)$ is isomorphic to  $(\mathbb{C}^\times)^k$ for some $k\in\mathbb{Z}_+$.
\end{co}

Recall that the subgroup generated by $\exp\mathfrak{h}$ for a Lie subalgebra $\mathfrak{h}$  of the Lie algebra~$\mathfrak{g}$ corresponding
to $G$ is said to be  \emph{integral}; and, by definition, the  \emph{solvable radical} is the integral subgroup corresponding to the solvable radical of~$\mathfrak{g}$ \cite[p.\,591, Definition 16.2.1]{HiNe}. The following assertion is used in the proof of
Theorem~\ref{linarizc}.

\begin{thm}\label{linrad0} \cite[Theorem~2.1]{ArAnF}
A connected complex Lie group is linear if and only if its solvable
radical is linear.
\end{thm}
Note that in 1943 Malcev published his version of Theorem~\ref{linrad0} for real Lie groups (the linearity of the
Levi subgroup is required additionally, but the analogous condition holds automatically in the complex
case) in~\cite{Ma} (for example, see a proof in \cite[p.\,221,  Chapter XVIII, Theorem~4.2]{Ho65}). However, to the
author's knowledge, the complex case has not attracted the attention of specialists over the past decades
and was only studied in~\cite{ArAnF}.

First, we prove some consequence of Theorem~\ref{linrad0}.
\begin{thm}\label{linrad}
Let $G$ be a connected complex Lie group, and let $R$ be the solvable radical of $G$. Then
 $\LinC(R)=\LinC(G)$.
\end{thm}

\begin{proof}
The inclusion  $\LinC(R)\subset \LinC(G)$ is obvious. To prove the inverse inclusion we recall that
$\LinC(G)$ is central. Hence, $\LinC(R)$ is central too, and so it is normal in $G$. We claim that $G/\LinC(R)$ is linear.

Indeed, denote by $\mathfrak{g}$, $\mathfrak{r}$ and $\mathfrak{l}$ the Lie algebras associated with $G$, $R$, and $\LinC(R)$.
Since~$\mathfrak{l}$ is solvable,
the radical of $\mathfrak{g}/\mathfrak{l}$  coincides with $\mathfrak{r}/\mathfrak{l}$ (for example, see \cite[Lemma 4.10]{ArOld}). Since the radical of a Lie group is
the integral subgroup corresponding to the radical of the Lie algebra; therefore, as the exponential mapping
is functorial,  $R/\LinC(R)$ is the integral subgroup generated by $\mathfrak{r}/\mathfrak{l}$, whence $R/\LinC(R)$  coincides with the radical of  $G/\LinC(R)$.

The quotient group by the linearizer of a group is always linear   \cite[p.\,578, Lemma 15.2.14]{HiNe}; in
particular, $R/\LinC(R)$  is linear. Since the linearity of the radical of a connected complex Lie group
implies the linearity of the initial group by Theorem~\ref{linrad0}, we infer that $G/\LinC(R)$ is linear. Hence,
$\LinC(G)\subset\LinC(R)$.
\end{proof}

Furthermore, we need two auxiliary assertions.

\begin{lm}\label{Linqu}
Let $G$ be a complex Lie group, and let $H$ be a closed subgroup of $\LinC(G)$. Then
$\LinC(G/H)=\LinC(G)/H$.
\end{lm}
\begin{proof}
Note that $\LinC(G/H)\subset \LinC(G)/H$ because $$(G/H)/(\LinC(G)/H)\cong G/\LinC(G)$$ and the latter group is linear.

On the other hand, assume that  $g\in \LinC(G)$. If $\alpha$ is a finite-dimensional holomorphic representation
of  $G/H$, then the composition of $\alpha$ with $G\to G/H$  is a finite-dimensional holomorphic representation
of $G$ which sends $g$ to $1$. Thus, the kernel of $\alpha$  contains $gH$. Since $\LinC(G/H)$ is the intersection of the
kernels of all such representations, $gH\in \LinC(G/H)$.
Hence, $\LinC(G)/H \subset \LinC(G/H)$.
\end{proof}

\begin{pr}\label{Linconn}
If a complex Lie group $G$ is connected then so is $\LinC(G)$.
\end{pr}
\begin{proof}
Let $\Lambda_0$ be a connected component of $\LinC(G)$. Putting $H=\Lambda_0$ in Lemma~\ref{Linqu}, we see that  $\LinC(G/\Lambda_0)\cong \LinC(G)/\Lambda_0$. Since $\LinC(G)/\Lambda_0$  is discrete, it suffices to prove that the linearizer of
a connected Lie group is discrete if and only if it is trivial.

So, let  $\LinC(G)$  be discrete. Recall that every solvable subgroup of a connected locally compact
group is compactly generated (see \cite{HN09}).  In particular, this holds for abelian subgroups, and so $\LinC(G)$ is finitely generated since it is discrete.  As it is easy to see, every finitely generated abelian group is
residually finite. In particular, given a nonunit  $g\in \LinC(G)$, there is a subgroup $\Lambda_1$ of $\LinC(G)$ such that $\LinC(G)/\Lambda_1$  is finite and  $g\notin\Lambda_1$. Moreover, we may assume that  $\LinC(G)/\Lambda_1$  is a finite cyclic group.
Passing to $G/\Lambda_1$ and applying Lemma~\ref{Linqu}, we reduce the claim to the case when $\LinC(G)$  is finite.

Assume now that  $\LinC(G)$  is finite.  In this case, since $G/\LinC(G)$  is connected and linear, and
 $\LinC(G)$ is central, $G$ is linear \cite[Corollary 16.3.9]{HiNe}; i.e., $\LinC(G)$  is trivial.
\end{proof}

\begin{rem}
The linearizer of a connected real Lie group (in the sense of real Lie groups) need not
be connected. For example, it is well known that the linearizer of the universal covering group of $\mathrm{SL}(2,\mathbb{R})$ is isomorphic to~$\mathbb{Z}$.
\end{rem}

\begin{proof}[Proof of Theorem~\ref{linarizc}]
Notice first that the Morimoto subgroups of $G$ and $\LinC(G)$ coincide.
Indeed, denote them by $M$ and $M_1$, respectively. Since every holomorphic function on $M$ is constant, so
are all holomorphic representations. Hence, $M\subset \LinC(G)$.   Since $G/M$ is a Stein group, so is its every Lie
subgroup; in particular,  $\LinC(G/M)$. Putting $H = M$ in Lemma~\ref{Linqu}, we get $\LinC(G/M)=\LinC(G)/M$. Hence,  $\LinC(G)/M$  is a Stein group. Using the definition of $M_1$ we infer that $M_1\subset M$. On the other
hand, it follows from the explicit construction of the Morimoto subgroup \cite[\S\,2]{Mo65}  that $M_1$ is maximal
among those Lie subgroups of $\LinC(G)$  that have only constant holomorphic functions on $\LinC(G)$. Thus, $M_1=M$.

By Theorem~\ref{linrad}, we may assume that $G$ is solvable. By \cite[Theorem 3.2]{Mo66} $\LinC(G)$  is of the form~\eqref{Mdec}, and $M$ coincides with the Morimoto subgroup of $G$, as was just shown. It remains to demonstrate that
$m=0$.  Assume the contrary. Then  $\LinC(G)\cong N\times \mathbb{C}$, where $N\cong M \times (\mathbb{C}^\times)^k\times \mathbb{C}^{m-1}$. Since $N$ is central,
we may consider the quotient group $G/N$. Putting $H = N$ in Lemma~\ref{Linqu}, we get $\LinC(G/N)\cong \LinC(G)/N$. So, it suffices to consider the case when $\LinC(G)$ is isomorphic to~$\mathbb{C}$.

Assume that $\LinC(G)\cong\mathbb{C}$.  Denote $G/\LinC(G)$ by $G_1$ and the quotient homomorphism $G\to G_1$ by~$\sigma$. Let $K$ be a maximal compact subgroup of $G$. Put $K_1=\sigma(K)$.  Then by  \cite[p.\,545, Theorem 14.3.13]{HiNe} $K_1$ is a maximal compact subgroup of $G_1$. Furthermore, it follows from  \cite[p.\,540, Lemma 14.3.3(4)]{HiNe}  that
a maximal compact subgroup of a solvable Lie group is abelian and, consequently, it is a torus. Thus,
we have a homomorphism  $\sigma'\!:K\to K_1$ of tori.

From the exactness of the homotopy sequence of pointed bundle (for example, see \cite[Lecture~5, pp.\,219--220, Formula~(15) and Proposition~2]{Po84}) and from the contractibility of $\LinC(G)$ we infer that $\sigma$ induces an isomorphism   $\pi_1(G)\cong \pi_1(G_1)$ of the fundamental groups.  It follows from the second splitting theorem \cite[p.\,544, Theorem 14.3.11]{HiNe}  that $G$ and $G_1$ are homotopically equivalent to $K$ and $K_1$. Hence,  $\pi_1(K)\cong \pi_1(K_1)$.  Since  $K\cong \mathbb{T}^n$ and $K_1\cong \mathbb{T}^{n_1}$ for some $n,n_1\in\mathbb{Z}$,  $\pi_1(K)\cong \mathbb{Z}^n$ and $\pi_1(K_1)\cong \mathbb{Z}^{n_1}$. It is easy to verify
that every homomorphism $\mathbb{T}^n\to \mathbb{T}^{n_1}$ is uniquely defined by the corresponding homomorphism $\mathbb{Z}^n\to\mathbb{Z}^{n_1}$. In particular, this implies that $\sigma'$  is an isomorphism.

Note that $G$ and $G_1$ are Stein groups. Indeed, the Morimoto subgroup lies in the linearizer. In our
case, it is either trivial or isomorphic to~$\mathbb{C}$. By the connectedness of Morimoto subgroups \cite[Lemma~4]{Mo65}, we infer that they are trivial for both groups, as required.

Denote by $K^*$ and $K_1^*$ the smallest complex integral subgroups containing $K$ and $K_1$, respectively.
These subgroups are closed \cite[Corollary 14.5.6]{HiNe}; and since $G$ and $G_1$ are Stein groups, $K^*$ and $K_1^*$ have
no compact factors. Consequently, these are the universal complexifications of  $K$ and $K_1$  \cite[Proposition 15.3.4 (i)]{HiNe}. Since universal complexification is a functor, $\sigma$ maps $K^*$ onto $K_1^*$ isomorphically.

Since $G/\LinC(G)$  is connected and linear, there is a decomposition into a semidirect product:
$G/\LinC(G) = B_1\rtimes K_1^*$, where $B_1$  is simply connected and solvable \cite[p.\,601, Theorem 16.3.7]{HiNe}. Put $B\!:=\sigma^{-1}(B_1)$. Obviously, $\LinC(G)\subset B$, and $\LinC(G)$ coincides with the kernel of the restriction of $\sigma$ to $B$. Thus, $B/\LinC(G)\cong B_1$; and, consequently, $B$ is also simply connected and solvable (since $\LinC(G)\cong\mathbb{C}$).

We claim that $G=B\rtimes K^*$. Indeed, firstly $\sigma(g)=b_1 l_1$ for every $g\in G$, where $b_1\in B_1$ and $l_1\in K_1^*$. Let $l\in K^*$ be the preimage of $l_1$. Then $g=(gl^{-1})l$, and $gl^{-1}\in B$ since $\sigma(gl^{-1})=b_1$. Thus $G=BK^*$.  Secondly, if  $g\in B \cap K^*$, then $\sigma(g)\in B_1 \cap K_1^*=\{1\}$,  whence $g=1$ since $K^*\to K_1^*$ is an isomorphism. Thus,  $B \cap K^*=\{1\}$. So $G=B\rtimes K^*$ with $B$ simply connected and solvable, and $K$ is
a maximal compact subgroup. Applying in the inverse direction Theorem 16.3.7 from \cite{HiNe}, we infer that $G$
is linear. The so-obtained contradiction with the nontriviality of the linearizer completes the proof of
Theorem~\ref{linarizc}.
\end{proof}

\subsection*{Acknowledgment} The author thanks the referee for the useful comments that improved the
presentation.

\end{document}